\documentclass[11pt]{amsart}
\usepackage{epsf}

\newcommand{\A}{{\mathfrak A}}

\newcommand{\gl}{{\mathfrak{gl}}}
\newcommand{\wS}{\widetilde S}
\newcommand{\wR}{\widetilde R}
\newcommand{\wPsi}{\widetilde\Psi}
\newcommand{\wa}{\widetilde a}
\newcommand{\wxi}{\widetilde\xi}
\def\O{\mathcal O}

\def\matho#1{\mathop{\mathrm{#1}}}

\newcommand*{\ato}[2]{{\genfrac{}{}{0pt}{}{#1}{#2}}}
\def\suml{\sum\limits}
\DeclareMathSizes{11.1}{10}{8}{6}

\newcommand\nfrac[2]
{\dfrac{\raisebox{-1pt}{\fontsize{11.1}{10pt}\selectfont$#1$}}
       {\raisebox{ 1pt}{\fontsize{11.1}{10pt}\selectfont$#2$}}}

\newcommand\ii[1]{$\mathrm{(iii)}_{#1}$}

\newcommand{\Alt}{\matho{Alt}}

\newcommand{\Altl}{\Alt\limits}
\newcommand{\Tr}{\matho{Tr}}

\newcommand{\ad}{\matho{ad}}
\newcommand{\fin}{\mathrm{fin}}

\newcommand{\even}{\mathrm{even}}
\newcommand{\inner}{\mathrm{inner}}
\newcommand{\Circle}{\mathrm{Circle}}
\newcommand{\Dif}{\mathrm{Dif}}
\newcommand{\Poiss}{\matho{Poiss}}

\newcommand{\C}{\mathbb C}

\newtheorem*{theorem}{Theorem}
\newtheorem*{lemma}{Lemma}

\newtheorem*{conjecture}{Conjecture}

\newtheorem*{fprobl}{First Problem}
\newtheorem*{sprobl}{Second Problem}
\newtheorem*{klemma}{Key Lemma}

\theoremstyle{remark}
\newtheorem*{remark}{Remark}
\newtheorem*{example}{Example}

\theoremstyle{definition}
\newtheorem*{defin}{Definition}

\author{Boris Shoikhet}
\title[Lifting formulas II]%
{Lifting
formulas II}
\date{1998}
\address{IUM, 11 Bol'shoj Vlas'evskij per.,
Moscow 121002, Russia}
\email{borya@mccme.ru}

\begin{document}
\maketitle
 \sloppy

\def\nfracp#1#2{\nfrac{\partial#1}{\partial#2}}

\section*{Introduction}

Let $\A$ be an associative algebra, we will also think of~$\A$
as a Lie algebra with the bracket $[a,b]=a\cdot b-b\cdot a$.
Suppose that there is a trace $\Tr\colon\A\to\C$ on the
associative algebra~$\A$ (this means that $\Tr([A,B])=0$ for
any $A,B\in\A$) and a set of its (exterior) derivations
$\{D_1,\dots,D_n\}$; such that the following conditions 1)--2) hold:

1) $\Tr(D_iA)=0$ for all $i=1,\dots,n$ and all $A\in\A$;

2) $[D_i,D_j]=0$ for all $i$, $j$.

It was proven in \cite{Sh} (and probably elsewhere as well)
that if conditions 1) and 2) hold, then
\begin{equation}
\Psi_{n+1}(A_1,\dots,A_{n+1})=\Altl_A\Altl_D\Tr
(D_1A_1\cdot\ldots\cdot D_nA_n\cdot A_{n+1})
\end{equation}
is a $(n+1)$-cocycle on the Lie algebra~$\A$.

Furthermore, let us suppose that instead of condition~2) we have
the weaker condition $\text{(i)}+\text{(ii)}$:

(i) $[D_i,D_j]=\ad Q_{ij}$  for all $i,j=1,\dots,n$\ \
($Q_{ij}\in\A$ and $Q_{ij}=-Q_{ji}$)

(ii) $\Altl_{i,j,k}D_k(Q_{ij})=0$ for all $i$, $j$, $k$.

Condition (ii) is analogous to the Jacobi identity. In fact, one
can deduce from the Jacobi identity that
$\Alt_{i,j,k}D_k(Q_{ij})$ lies in the center of the Lie
algebra~$\A$.

When conditions 1), (i) and (ii) hold,
a
formula for the $(n+1)$-cocycle on the Lie algebra~$\A$
was found
in \cite{Sh};
this formula is
quantization on~$Q_{ij}$ of the formula~(1). This quantization
is called the \emph{lifting formula}. The fact that the lifting
formula in reality defines a cocycle was  the  Main  Conjecture
in~\cite{Sh}.

The present paper contains the solution of the following two problems:

\begin{fprobl}
When   conditions  \emph{1)--2)}  above  hold  $([D_i,D_j]=0)$,  to
construct formulas for $n+3$, $n+5$, $n+7$, \dots ---  cocycles,
which are analogous to the formula~\emph{(1)}.
\end{fprobl}

\begin{sprobl}
To construct a quantization of these formulas
on~$Q_{ij}$, when condition~\emph{2)} is replaced by the
condition $(i)+(ii)$.
\end{sprobl}

The First Problem is solved completely in Section~1; the answer
to the Second Problem is given in Theorem~2.3; we prove this
Theorem as well as the Main Conjecture from~\cite{Sh} in Section~3.

Conditions $(1)+(i)+(ii)$ arise in the
following situation. Let $\A=\Psi\Dif_n(S^1)$ be the associative
algebra of formal pseudodifferential operators on $(S^1)^n$
(see~\cite{A}). It is easy to check (see also~\cite{A}) that the
bracket of two such operators has coefficient zero at the term
$x_1^{-1}\cdot\ldots\cdot
x_n^{-1}\cdot\partial_1^{-1}\cdot\ldots\cdot\partial_n^{-1}$ (in
any coordinate system). Therefore this coefficient defines a
trace functional $\Tr\colon\Psi\Dif_n(S^1)\to\C$ which is known
as a ``noncommutative residue.'' Furthermore, $\ad(\ln x_1)$,
\dots, $\ad(\ln x_n)$, $\ad(\ln\partial_1)$, \dots,
$\ad(\ln\partial_n)$ gives us a set of $2n$ exterior derivations
on the associative algebra $\Psi\Dif_n(S^1)$ satisfying
conditions~1), (i) and~(ii). It was shown in \cite{Sh} that
$$
[\ln\partial,\ln x]=\ad\left(x^{-1}\partial^{-1}+\nfrac12x^{-2}\partial^{-2}+
\nfrac23x^{-3}\partial^{-3}+\ldots+\nfrac{(n-1)!}nx^{-n}\partial^{-n}+\ldots
\right)
$$

The $(2n+1)$-cocycle constructed in~\cite{Sh} is not
cohomologous to zero; it remains non\-cohomologous to zero when
restricted to the Lie subalgebra $\Dif_n\hookrightarrow\Psi\Dif_n(S^1)$
of the (polynomial) differential operators on~$\C^n$. The
simplest way to check this fact is to construct its limit
with respect to the deformation of the Lie algebra $\Poiss_{2n}$
and to prove that this $(2n+1)$-cocycle is not cohomologous to zero.

The characteristic feature of our situation is that the $2n+3$,
$2n+5$, $2n+7$, \dots --- cocycles on the Lie algebra $\Dif_n$
constructed in this paper have quite complicated Hamiltonian
limits; the question about their cohomological nontriviality
remains  open.

One final comment: it was proved in \cite{FT} that
$$
H^*(\gl_\infty^\fin(\Dif_n);\C)=\wedge^*(\xi_{2n+1},\xi_{2n+3},\xi_{2n+5},
\dots);
$$
furthermore, the trace on the associative  algebra $\Psi\Dif_n(S^1)$ can
be extended to the associative  algebra
$\gl_\infty^\fin(\Psi\Dif_n(S^1))$ and the derivations $\{\,\ln
x_i,\ln\partial_i,\,i=1,\dots,n\,\}$  also act on this algebra.
This fact allows us to construct the corresponding cocycles
on the Lie algebra $\gl^\fin_\infty(\Dif_n)$, and our cocycles are
pullbacks of these cocycles
with
respect to the inclusion
$\Dif_n\hookrightarrow\gl_\infty^\fin(\Dif_n)$\ \ ($D\mapsto
D\cdot E_{11}$). This construction together with the above
result from~\cite{FT} imply that:

(1) Automorphisms of the Lie algebra $\Psi\Dif_n(S^1)$ act
trivially on our cocycles;

(2) Another choice of the exterior derivations does not lead
to new cohomological classes.

Actually, our methods allow us to construct just one unique cohomology
class
of the Lie algebra $\Dif_n$ in each odd dimension $\ge2n+1$.
\smallskip

I am grateful to B.L.\,Feigin for numerous valuable discussions.

\section{Case $[D_i,D_j]=0$}

\subsection{}
Let $\A$ be an associative  algebra with trace functional
$\Tr\colon\A\to\C$, and let $\{D_1,\dots,D_n\}$ be a set of its
(exterior) derivations  satisfying the following condition:
\begin{equation}
\Tr(D_iA)=0\quad\text{for any $i=1,\dots,n$ and any  $A\in\A$}.
\end{equation}

Let $l\ge1$ be an integer and let us consider the following expression:
\begin{equation}
S_{a_1,\dots,a_{n+2l}}(A_1,\dots,A_{n+2l})=\Altl_{A,D}\Tr(T_1
\cdot\ldots\cdot T_{n+2l})
\end{equation}
where

(i) $A_1,\dots,A_{n+2l}\in\A$;

(ii) $a_1,\dots,a_{n+2l}\in\{0,1\}$;

(iii) $T_i=D_{j(i)}A_i$ if $a_i=1$\ \ ($j(i)$ is defined below)
$T_i=A_i$ if $a_i=0$;

(iv) $a_1=1$ and $T_1=D_1A_1$;

(v) $j$ takes the values from $1$ to $n$ \emph{in turn}, i.e., if
$i_1<i_2$ then $j(i_1)<j(i_2)$.

\begin{defin}
\begin{equation}
S_\even=\sum_{(a_1,\dots,a_{n+2l})\in a_\even} S_{a_1,\dots,a_{n+2l}}
(A_1,\dots,A_{n+2l})
\end{equation}
where $a_\even$ is the set of the sequences with values in
$\{0,1\}$ of length $n+2l$ and such that

(i) $a_1=1$;

(ii) $n$ of the $a_i$'s are equal to~$1$ and $2l$ of the $a_i$'s are
equal to~$0$.

(iii) if $a_i=1$, $a_j=1$, $j>i$ and  $a_{i+1},\dots,a_{j-1}=0$,
then  $j-i$  is  \emph{odd}.  In  other  words,  there are an
\emph{even} number of $0$'s between the two nearest~$1$'s.
\end{defin}

\begin{remark}
The condition (iii) should also hold  for the ``\emph{tail}''  of
the  sequence  $(a_j)$,  as if the numbers $a_i$'s were placed on a
circle.
\end{remark}

\begin{lemma}
If $[D_i,D_j]=0$ for all $i,j=1,\dots,n$ then $S_\even\equiv0$.
\end{lemma}

\subsection{}
\begin{proof}[Proof of Lemma \emph{1.1}]

Let us assign to  each  expression  $S_{a_1,\dots,a_{n+2l}}$,
$(a_i)\in  a_\even$  an  expression  $\wS_{a_1,\dots,a_{n+2l}}$
($A_1,\dots,A_{n+2l}$) in the following way:
\begin{equation}
\wS_{a_1,\dots,a_{n+2l}}=\sum_{j:a_j=1}\Altl_{A,D}
\Tr(D_j(T_{1,j}\cdot\ldots\cdot T_{n+2l,j}))
\end{equation}
where $T_{i,j}=T_i$ for $i\ne j$ and $T_{j,j}=A_j$. Then
Lemma~1.1 follows from the following statement:

\begin{lemma}
$\suml_{(a_i)\in a_\even}\wS_{a_1,\dots,a_{n+2l}}=(n+l)\cdot S_\even$.
\end{lemma}

Lemma 1.1 follows from this Lemma because of identity~(2).

\begin{proof}[Proof of Lemma]
Each summand in $\wS_{a_1,\dots,a_{n+2l}}$ can be written as a further sum
by the Leibniz rule. We have:
\begin{equation}
\Altl_{A,D}D_j(T_{1,j}\cdot\ldots\cdot T_{n+2l,j})=
\sum_{i:T_{i,j}=A_i}\Altl_{A,D}\Tr(T_{1,j}\cdot\ldots
\cdot D_jT_{i,j}\cdot\ldots\cdot T_{n+2l,j})
\end{equation}
because $[D_i,D_j]=0$ for all~$i$. The summand in the r.h.s.\
of~(6) corresponding to $i=j$ is equal to
$S_{a_1,\dots,a_{n+2l}}$. As $j$ ranges from~$1$ to~$n$ we
obtain $n\cdot S_{a_1,\dots,a_{n+2l}}$. The other summands in
the r.h.s.\ of~(6) are of the following three types:

(i) Summands of the form $\Altl_{A,D}\Tr(T_{1,j}\cdot\ldots\cdot
D_jT_{i,j}\cdot\ldots T_{n+2l,j})$\ \ ($T_{i,j}=A_i$) for which,
when $i<j$, there exists~$k$ such that $i<k<j$; when $i>j$, there
exists $k$ such that $j<k<i$;
and such that $a_k=1$.

\begin{figure}[h]
\centerline{\epsfbox{pict.1}}
\caption{}
\end{figure}

We move $1$ from the $j$-th place to the $i$-th place, where it
divides the sequence of
zeroes of  even length into two sequences. One of these
sequences is necessarily of odd length, while the second
sequence has even length. In Fig.~1, $1$'s are located on black
squares before and after moving, and $0$'s are on white squares.

For each summand of this type there is a ``dual'' summand;
the sum of a summand and its dual is~$0$. For example,
the dual of Fig.~1 is Fig.~2 below:

\begin{figure}[h]
\centerline{\epsfbox{pict.2}}
\caption{}
\end{figure}

(ii) However, the dual summand for a summand of type~(i)
is \emph{not} necessarily also of type~(i). The typical situation
is shown in Fig.~3:

\newcounter{tmp}
\setcounter{tmp}{\value{equation}}
\setcounter{equation}{0}
\renewcommand{\theequation}{\roman{equation}}
\begin{figure}[h]
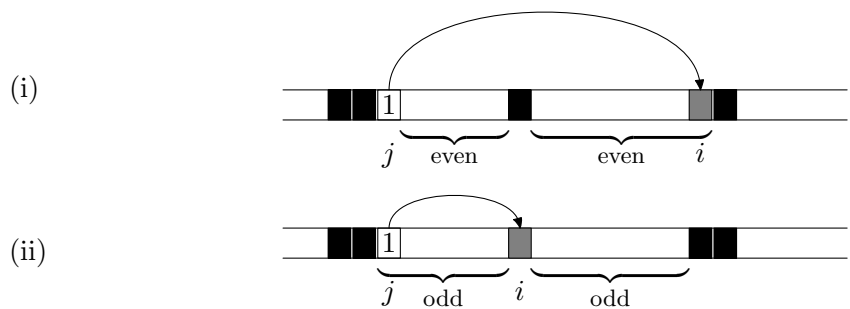

\begin{equation}
\vcenter{\epsfbox{pict.3}}
\end{equation}
\begin{equation}
\vcenter{\epsfbox{pict.4}}
\end{equation}
\caption{The summand of the type~(i) and the dual summand of the type~(ii)}
\end{figure}
\setcounter{equation}{\value{tmp}}
\renewcommand{\theequation}{\arabic{equation}}

Figure~3 correspons to the case when there
does \emph{not} exist $a_k=1$ between
$a_i=0$ and where $a_j=0$ and $a_j$ divides the corresponding sequence
into two odd parts. The two dual summands cancel one another
because of the alternation of the~$D$'s.

(iii) There remain summands of a third type, which is
shown schematically in Fig.~4. In this case there does not exist
$a_k=1$ between $a_i=0$ and $a_j=0$, and $a_j$ divides the
corresponding even sequence into two \emph{even} parts. We
differentiate types \ii1 and \ii2; strictly speaking, type \ii1 is
a special case of type \ii2.

\setcounter{tmp}{\value{equation}}
\setcounter{equation}{0}
\makeatletter
\def\@eqnnum{{\normalfont \normalcolor $(iii)_{\arabic{equation}}$}}
\makeatother
\begin{figure}[h]
\begin{equation}
\vcenter{\epsfbox{pict.5}}
\end{equation}
\begin{equation}
\vcenter{\epsfbox{pict.6}}
\end{equation}
\caption{Types \ii1 and \ii2}
\end{figure}
\setcounter{equation}{\value{tmp}}
\makeatletter
\def\@eqnnum{{\normalfont \normalcolor (\theequation)}}
\makeatother

Summands of type (iii) do \emph{not} have dual summands, but
they are of the form $S_{a_1,\dots,a_{n+2l}}$ for some $(a_i)\in
a_\even$, therefore they are summands of~$S_\even$. These
summands together with the summands of the form
$\Altl_{A,D}\Tr(T_{1,j}\cdot\ldots\cdot
D_jT_{j,j}\cdot\ldots\cdot T_{n+2l,j})$ lead to the
coefficient $(n+l)$ in the statement of the Lemma.
\end{proof}
\let\qed\relax
\end{proof}

\subsection{}
Our aim is to construct the $(n+2l-1)$-cocycle on the Lie
algebra~$\A$ analogous to the $(n+1)$-cocycle given by
formula~(1).

Let us define the $(n+2l-1)$-cochain
$R_{a_1,\dots,a_{n+2l}}(A_1,\dots,A_{n+2l-1})$ for each
$S_{a_1,\dots,a_{n+2l}}$, $(a_i)\in a_\even$. Roughly
speaking, we just shorten any sequence of consecutive zeroes to a single zero;
after this procedure, the sequence will have
odd length. We choose the first such sequence. More
precisely, let us define the sequence
$(\wa_i)_{i=1,\dots,n+2l-1}$ for every sequence $(a_i)\in
a_\even$ in the following way:

let $s_1=\min\limits_i(a_i=0)$, $s_2=\max\limits_{i>s_1}(a_i=1)$; then:

$\wa_1,\dots,\wa_{s_1-1}=1$;

$\wa_{s_1},\dots,\wa_{s_2-2}=0$;

$\wa_i=a_{i+1}$ for $s_2-1\le i\le n+2l-1$.

Let $R_{a_1,\dots,a_{n+2l}}(A_1,\dots,A_{n+2l-1})=
\Altl_{A,D}\Tr(P_1\cdot\ldots\cdot P_{n+2l-1})$, where

$P_i=D_{j(i)}A_i$ for $\wa_i=1$,

$P_i=A_i$ for $\wa_i=0$

and $j(i)$ is defined in the following way: $j(1)=1$;
$j(i_1)<j(i_2)$ when $i_1<i_2$ and $j=1,\dots,n$. In other words,
$j$ takes values from $1$ to $n$ \emph{in turn}.

\begin{lemma}
$$
d(R_{a_1,\dots,a_{n+2l-1}})(A_{n+2l},A_1,A_2,\dots,A_{n+2l-1})=
(-1)^{n+2l-s_1+1}\cdot S_{a_1,\dots,a_{n+2l}}(A_1,\dots,A_{n+2l})
$$
\end{lemma}

\begin{proof}
A direct calculation in the spirit of~\cite{Sh}, Sections 1--%
3.
\end{proof}

\begin{theorem}
\begin{equation}
\Psi_{n,l}^0=\sum_{(a_i)\in a_\even}(-1)^{s_1}R_{a_1,\dots,a_{n+2l}}
(A_1,\dots,A_{n+2l-1})
\end{equation}
is a $(n+2l-1)$-cocycle on the Lie algebra~$\A$.
\end{theorem}

\begin{proof} Follows from Lemma 1.1.
\end{proof}

\section{Quantization}

\subsection{}
Let $\A$ be the associative algebra with trace~$\Tr$;
$D_1,\dots,D_n$ are its (exterior) derivations, which satisfy the
following conditions:
\begin{gather}
[D_i,D_j]=\ad Q_{ij};\quad Q_{ij}=-Q_{ji}\\
\Altl_{i,j,k}D_k(Q_{ij})=0
\end{gather}
In the case when $Q_{ij}=0$ for all $i,j$ we have constructed in Section~1
the $(n+2l-1)$-cocycle on the Lie algebra~$\A$, for
all integral $l\ge1$. Our next problem is to quantize
this cocycle and to find
the $(n+2l-1)$-cocycle for all $Q_{ij}\in\A$
satisfying conditions~(8) and~(9).

In \cite{Sh} such a cocycle was constructed for $l=1$
(quantization of formula~(1)). Let us recall that
construction.

Suppose that $n\ge2$. We will consider intervals of length
$n-2$ with some \emph{marked} integral points, such that the
\emph{distance between any two marked points greater or equal
than~$2$}. Let us denote by $I_k$ the set of all such intervals
with $k$ marked points ($1\le k\le\left[\nfrac n2\right]$).
Denote by $1,\dots, n-1$ the integral points of the interval.

\begin{defin}
Suppose that $t\in I_k$ and $i_1<\ldots<i_k$ are its marked points
($1\le i_1,i_k\le\left[\nfrac n2\right]$ and $i_{s+1}-i_s\ge2$ for all
$s=1,\dots,k-1$). Then
\begin{equation}
\O(t)=\Altl_{A,D}\Tr(P_{1,t}\cdot\ldots\cdot P_{n+1,t})
\end{equation}
where

{\multlinegap=0pt
\begin{multline*}
\begin{split}
&\left.
\begin{split}
&P_{i,t}=A_j\cdot Q_{j,j+1}\\
&P_{j+1,t}=A_{j+1}
\end{split}
\right\}\
\text{when $j=i_s$ for some $s=1,\dots,k$ (i.e., if point~$j$ is
marked)}\\
&\,P_{i,t}=D_jA_j\ \text{when $j$ and $j-1$ are not marked and $j\ne n+1$}
\\
&\,P_{n+1,t}=A_{n+1}
\end{split}\\
\end{multline*}}
\end{defin}

\begin{example}
If $n=6$ and $t=\vcenter{\epsfbox{pict.7}}$, $t\in I_2$, then
$$
\O(t)=\Altl_{A,D}\Tr(A_1\cdot Q_{12}\cdot A_2\cdot
D_3A_3\cdot A_4\cdot A_5\cdot D_6A_6\cdot A_7).
$$
If $t=\vcenter{\epsfbox{pict.8}}$, $t\in I_3$, then
$$
\O(t)=\Altl_{A,D}\Tr(A_1\cdot Q_{12}\cdot A_2\cdot
A_3\cdot Q_{34}\cdot A_4\cdot A_5\cdot Q_{56}\cdot A_6\cdot A_7).
$$
\end{example}

\begin{remark}
In formula (10) we are using the convention that
$Q_{ij}$\,``$=$''\,$[D_i,D_j]$ in the sense of the alternation
on~$D$. Furthermore, according to our convention we \emph{%
don't alternate the symbols $i$ and~$j$ in~$Q_{ij}$}.
\end{remark}

\begin{theorem}[Main Conjecture from \cite{Sh}]
Let $\Sigma_k=\suml_{t\in I_k}\O(t)$.

If conditions \emph{(8)} and \emph{(9)} hold, then
\begin{equation}
\Psi_{n,1}=\Altl_{A,D}\Tr(D_1A_1\cdot D_2A_2\cdot\ldots\cdot
D_nA_n\cdot A_{n+1})+\Sigma_1+\Sigma_2+\ldots+\Sigma_{\left[\nfrac n2\right]}
\end{equation}
is a $(n+1)$-cocycle on the Lie algebra~$\A$.
\end{theorem}

We will prove this Theorem in Section~3.

\subsection{}
Here we give a motivation for Theorem~2.1 and prove the
Key Lemma, which we will use in Subsection~2.3 in order to
formulate Theorem~2.3, which generalizes Theorem~2.1 to the case
of arbitrary $l\ge1$. We will prove this more general Theorem in
Section~3.

Let us suppose that all derivations~$D_i$ in the formula~(7)
for a $(n+2l-1)$-cocycle on the Lie algebra~$\A$ are \emph{inner},
i.e., $D_iA=D_i\cdot A-A\cdot D_i$,\ \ $D_i\in\A$; now we do not
suppose that the condition $[D_i,D_j]=0$ hold. We define in this
case the $(n+2l-1)$-\emph{cochain}~$\Psi^\inner_{n,l}$ by the
formula~(7). In the same way, we define the
$(n+2l)$-\emph{cochain}~$S_\even^\inner$ by the formulas~(3),~(4).
It is clear that
\begin{equation}
d\Psi_{n,l}^\inner=S_\even^\inner.
\end{equation}
Now we replace all factors of the form $D_iA_j$ by
$D_i\cdot A_j-A_j\cdot D_i$ in the formula and remove all parentheses. We
obtain:
\begin{equation}
\Psi^\inner_{n,l}=\wPsi^\inner_{n,l}+r_{n,l}
\end{equation}
where $\wPsi_{n,l}^\inner$ is the sum of all summands
in~$\Psi_{n,l}^\inner$ in which no two $D_i$ and $D_j$ are
consecutive, and $r_{n,l}$ is the sum
of the remaining summands in~$\Psi_{n,l}^\inner$, i.e., those of the form
$\ldots\cdot D_i\cdot D_j\cdot\ldots$
for some~$i,j$.

\begin{klemma}
$\wPsi^\inner_{n,l}$ is a $(n+2l-1)$-cocycle on the Lie
algebra~$\A$.
\end{klemma}

\begin{remark}
In fact, $\wPsi^\inner_{n,l}$ is a coboundary.
\end{remark}

\begin{proof}[Proof of Key Lemma]
This Lemma is similar to Lemma~1.1, however, formally it is
another statement. By analogy with~(13), we have a similar
decomposition of~$S_\even^\inner$:
\begin{equation}
S_\even^\inner=\wS_\even^\inner+r_\even.
\end{equation}

\begin{lemma}
\leavevmode\par
\emph{(i)} $d\wPsi_{n,l}^\inner=\wS_\even^\inner$;

\emph{(ii)} $dr_{n,l}=r_\even$.
\end{lemma}

\begin{proof}
Follows from formula~(12).
\end{proof}

Now all that remains is to prove that
$\wS^\inner_\even\equiv0$ (without any conditions about the
\emph{inner} derivations~$D_i$). This is a direct calculation in
the spirit of Lemma~1.1 (see Lemma~4.2 of~\cite{Sh} for the case
$l=1$).
\end{proof}

Furthermore, in the case when all $D_i$ are inner derivations
satisfying condition~(8), $r_{n,l}$ is represented by the
sum~$\Sigma$ of the terms of the form~$\O(t)$, where $t$ is
an interval with marked points (see~(10)). We increase $D_i\cdot
D_j$ by $[D_i,D_j]$ via the alternation. Hence
$\Psi^\inner_{n,l}-\Sigma$ is a cocycle (which, in fact, is
cohomologous to zero). If now the $D_i$'s are arbitrary
derivations
(not
necessary inner)
satisfying condition~(8), then it
turns out that this expression continues to be a cocycle,
\emph{whenever condition~\emph{(9)} holds}. This was proved
in~\cite{Sh} in the case $l=1$ and $n=4$,
in the first nontrivial case. In the case $l=1$
this reasoning leads to the formulation of Theorem~2.1.

\subsection{Case of arbitrary $l\ge1$}

We have:
\begin{equation}
\wPsi^\inner_{n,l}=\sum_{(a_i)\in a_\even}(-1)^{s_1}\cdot
\wR_{a_1,\dots,a_{n+2l}}(A_1,\dots,A_{n+2l-1})
\end{equation}
where:
\begin{equation}
\wR_{a_1,\dots,a_{n+2l}}=\widetilde{\Altl_{A,D}\Tr}(P_1\cdot\ldots\cdot
P_{n+2l})
\end{equation}
and

$P_i=[D_{j(i)},A_i]$ when $\wa_i=1$;

$P_i=A_i$ when $\wa_i=0$;

$j(1)=1$ and the indices $j$ takes the values from $1$ to $n$
\emph{in term}
 (formula~(7)).

The sign ``$\widetilde{\phantom{\Alt\Tr}}$'' in the r.h.s.\
of~(16) means the sum  on all the terms, in which no two $D_i$,
$D_j$ are consecutive.

We quantize each summand $R_{a_1,\dots,a_{n+2l}}$ separately.

\begin{defin}
\leavevmode\par

(i) Denote by $\Circle_k^{a_1,\dots,a_{n+2l}}$ the set of all
circles with $n+2l-1$ integral points from which~$k$ ($1\le
k\le\left[\nfrac n2\right]$) are marked. The distance between any two marked
points is $\ge2$. The points are enumerated by $1,\dots,n+2l-1$.
Point $i$ may be marked only if $\wa_i=1$ and $\wa_{i+1}=1$ (we
suppose that $\wa_{n+2l}=\wa_1$).

(ii) For $t\in\Circle_k^{a_1,\dots,a_{n+2l}}$ we define
$\O(t)$ by analogy with Definition~2.1 after replacing the interval
with the circle.

(iii) 
$\Sigma_k^{a_1,\dots,a_{n+2l}}=\suml_{t\in\Circle_k^{a_1,\dots,a_{n+2l}}}
\O(t).$
\end{defin}

\begin{theorem}
Let $\A$ be an associative algebra with~$\Tr$, and let $D_1,\dots,D_n$ be its
\emph(exterior\emph) derivations, which satisfy
conditions~\emph{(8)}, and~\emph{(9)}. Then
$$
\Psi_{n,l}=\Psi_{n,l}^0+\sum_{(a_i)\in a_\even}(-1)^{s_1}\sum_{k\ge1}
\Sigma_k^{a_1,\dots,a_{n+2l}}
$$
is a $(n+2l-1)$-cocycle on the Lie algebra~$\A$.
\end{theorem}

\section{Proofs}

We prove here Theorem 2.1 and  Theorem~2.3.  Actually  we  prove
Theorem~2.1  only  (the  case~$l=1$) but one can easily generalize
the proof to the case of an arbitrary integer~$l\ge1$ (Theorem~2.3).

\subsection{}
Let $\A$ be an associative  algebra with trace functional $\Tr$,
and let $D_1,\dots,D_n$ be its derivations, \emph{which satisfy
conditions}~(8) and~(9).

According to Lemma 1.1, we have:
\begin{multline}
\begin{split}
\Altl_{A,D}\Tr&\{D_1(A_1\cdot D_2A_2\cdot\ldots\cdot D_nA_n\cdot
A_{n+1}\cdot A_{n+2})\\
&+D_2(D_1A_1\cdot A_2\cdot D_3A_3\cdot\ldots\cdot D_nA_n\cdot
A_{n+1}\cdot A_{n+2})\\
&+\hbox to 7cm{\,\dotfill\,}\\
&+D_n(D_1A_1\cdot\ldots\cdot D_{n-1}A_{n-1}\cdot A_n\cdot\cdot A_{n+1}\cdot 
A_{n+2})\}
\end{split}\\
=(n+1)\cdot\Altl_{A,D}\Tr(D_1A_1\cdot\ldots\cdot     D_nA_n\cdot
A_{n+1}\cdot A_{n+2})+\Phi_1
\end{multline}
where $\Phi_1$ consists of terms \emph{linear} on~$Q_{ij}$.
Let us denote
the l.h.s.\ of~(17)
by~$X_1$.
Our aim is to
represent~$\Phi_1$ \emph{modulo the trace of the full derivations} as a sum of
terms of the form $d(\O(t))$ (where $t$ is an interval
of length $2n-2$ with $1$ marked point(see Subsection~2.1) and
$d$ is the differential in the cochain complex of the Lie
algebra~$\A$) and the expression~$\Phi_2$, consisting of terms
\emph{quadratic in~$Q_{ij}$}, and so on.

We have:
\begin{multline}
\Phi_1=\Altl_{A,D}\Tr\{A_1\cdot Q_{12}\cdot A_2\cdot
D_3A_3\cdot\ldots-A_1\cdot A_2\cdot Q_{12}\cdot
D_3A_3\cdot\ldots\\
+A_1\cdot D_2A_2\cdot Q_{13}\cdot A_3\cdot
D_4A_4\cdot\ldots-A_1\cdot D_2A_2\cdot A_3\cdot Q_{13}\cdot
D_4A_4\cdot\ldots+\ldots-\ldots\}.
\end{multline}
Now suppose that
\begin{multline}
X_2=\Altl_{A,D}\Tr\{-D_3(A_1\cdot A_2\cdot Q_{12}\cdot A_3\cdot
D_4A_4\cdot\ldots)+D_2(A_1\cdot A_2\cdot Q_{13}\cdot A_3\cdot
D_4A_4\cdot\ldots)\\
-D_4(A_1\cdot D_2A_2\cdot A_3\cdot Q_{13}\cdot A_4\cdot
D_5A_5\cdot\ldots)+\ldots-\ldots\}.
\end{multline}
The first, second, and third summands in~(19) correspond to the
second, third, and fourth summands in~(18), respectively.

In general, the summands of~(18) are of the following two types:

(i) summands which contain products of the form $A\cdot Q\cdot A$
--- there are no summands in~(19) corresponding to the
summands of~(18) of this type;

(ii) the remaining summands in~(18) have the form
$$
\ldots\cdot D_iA_{i-1}\cdot D_iA_i\cdot Q_{j,i+1}\cdot
A_{i+1}\cdot D_{i+2}A_{i+2}\cdot\ldots
$$
or
$$
\ldots\cdot D_{i-1}A_{i-1}\cdot A_i\cdot Q_{j,i}\cdot D_{i+1}
A_{i+1}\cdot\ldots;
$$
the following summands of~(19)
$$
D_i(\ldots\cdot D_{i-1}A_{i-1}\cdot A_i\cdot Q_{j,i+1}\cdot
A_{i+1}\cdot D_{i+2}A_{i+2}\cdot\ldots)
$$
and
$$
D_{i+1}(\ldots\cdot D_{i-1}A_{i-1}\cdot A_i\cdot Q_{j,i}\cdot
A_{i+1}\cdot D_{i+2}A_{i+2}\cdot\ldots)
$$
correspond to the above two summands of~(18), respectively.

\begin{lemma}
\begin{equation}
X_2=2\Phi_1-d\left((2n+2)\cdot\sum_{t\in I_1}\O(t)\right)-\Phi_2,
\end{equation}
where $\Phi_2$ consists of expressions \emph{quadratic}
in~$Q_{ij}$ and $d$ is the differential in the cochain complex
of the Lie algebra~$\A$.
\end{lemma}

\begin{proof}
A straightforward calculation.
\end{proof}

In the sequel we will need an explicit form for the summands
of~$\Phi_2$: they may be of the following two types:
\begin{equation}
\begin{split}
(i)\hfil&\ldots\cdot DA\cdot\underline A\cdot DA\cdot\ldots\cdot
DA\cdot\underline{A\cdot   Q\cdot   A}\cdot   DA\cdot\ldots\cdot
DA\cdot\underline{A\cdot Q}\cdot DA\cdot\ldots\\
(ii)\hfil&\ldots\cdot DA\cdot\underline A\cdot DA\cdot\ldots\cdot
DA\cdot\underline{A\cdot   Q\cdot   A}\cdot   DA\cdot\ldots\cdot
DA\cdot\underline{Q\cdot A}\cdot DA\cdot\ldots
\end{split}
\end{equation}

\subsection{}
Starting from $\Phi_2$ we construct an expression~$X_3$ in the
following way: there are two different cases for expressions
given by the formulas~(21):

1) expressions which contain two factors of the form $A\cdot
Q\cdot A$. They have the following form:
$$
\ldots\cdot A\cdot Q\cdot A\cdot\ldots\cdot A\cdot Q\cdot A\cdot\ldots;
$$
there do not exist any corresponding summands in $X_3$;

2) to
$$
\xi=\ldots DA\cdot\underline{A_i}\cdot
DA\cdot\ldots\cdot DA\cdot\underline{A\cdot Q\cdot A}\cdot
DA\cdot\ldots\cdot DA\cdot A_k\cdot Q_{j,k}\cdot D_{k+1}A_{k+1}\cdot\ldots
$$
corresponds an expression
$$
\wxi=D_{k+1}(\ldots\cdot DA\cdot\underline{A_i}\cdot
DA\cdot\ldots\cdot DA\cdot\underline{A_k\cdot Q_{j,k}\cdot
A_{k+1}}\cdot DA\cdot\ldots)
$$
and to
$$
\xi=DA\cdot\underline{A_i}\cdot DA\cdot\ldots\cdot
DA\cdot\underline{A\cdot Q\cdot A}\cdot DA\cdot\ldots\cdot DA\cdot
D_{k-1}A_{k-1}\cdot\underline{Q_{jk}\cdot A_k}\cdot DA\cdot\ldots
$$
corresponds an expression
$$
\wxi=D_{k-1}(\ldots\cdot DA\cdot\underline{A_i}\cdot
DA\cdot\ldots\cdot DA\cdot\underline{A_{k-1}\cdot Q_{j,k}\cdot
A_k}\cdot\ldots)
$$

\begin{defin}
\begin{equation}
X_3=\sum_{\ato{\text{all summands of~$\Phi_2$}}{\text{of type 2}}}
\Altl_{A,D}\Tr\wxi
\end{equation}
\end{defin}

\begin{lemma}
$$
X_3=2\Phi_2-d\left((2n+2)\sum_{t\in I_2}\O(t)\right)-\Phi_3
$$
where $\Phi_3$ is a sum of expressions cubic in~$Q_{ij}$.\qed
\end{lemma}

\subsection{}

\begin{proof}[Proof of Theorem 2.1]
We define the sequence $X_1,X_2,X_3,\dots$ by formulas analogous to (20)
and~(22)
and prove Lemmas analogous to Lemma~3.1 and~3.2. All that remains
is to note
that every $X_i\equiv0$ by formula~(2).
\end{proof}

\begin{conjecture}
$$
H^*(\gl_\infty^\fin(\Dif_n);\C)=
\wedge^*(\Psi_{n,1},\Psi_{n,2},\Psi_{n,3},\dots)
$$
\end{conjecture}

\end{document}